\documentclass[12pt]{amsart}
\usepackage{amsthm}
\usepackage{amstext}
\usepackage{amssymb}
\usepackage{mathrsfs}
\usepackage{graphicx}
\usepackage{leftidx}
\usepackage[english]{babel}
\usepackage{enumerate}
\usepackage{xfrac}
\usepackage{mathtools}
\usepackage{xcolor}
\usepackage[all]{xy}
\usepackage{tikz}

\usepackage[T1]{fontenc}

\usepackage{mathpazo}
\usepackage{hyperref}

\allowdisplaybreaks

\theoremstyle{plain}
\newtheorem{theorem}{Theorem}[section]
\newtheorem*{theorem*}{Theorem}
\newtheorem{proposition}[theorem]{Proposition}
\newtheorem{corollary}[theorem]{Corollary}
\newtheorem{lemma}[theorem]{Lemma}

\theoremstyle{definition}
\newtheorem{definition}[theorem]{Definition}

\newtheorem{remark}[theorem]{Remark}

\theoremstyle{remark}

\numberwithin{equation}{section}

\newcommand{\N}{\mathbb N}

\newcommand{\R}{\mathbb R}

\newcommand{\dd}{\mathrm{d}}
\newcommand{\dist}{\mathbf{d}}
\newcommand{\Var}{\operatorname{Var}}
\newcommand{\V}{\operatorname{V}}
\newcommand{\vv}{\operatorname{v}}
\newcommand{\Par}{\operatorname{Par}}
\newcommand{\BV}{\mathsf{BV}}
\newcommand{\BVl}{\mathsf{BV}_{loc}}
\newcommand{\AC}{\mathsf{AC}}
\newcommand{\ACl}{\mathsf{AC}_{loc}}
\newcommand{\Csf}{\mathsf{C}}
\newcommand{\Dsf}{\mathsf{D}}
\newcommand{\Lsf}{\mathsf{L}}

\title{The speed measure and absolute continuity for curves in metric spaces}

\author{Sebastian Boldt}
\address{Fakultät f\"ur Mathematik\\Technische Universit\"at Chemnitz\\Germany}
\email{sebastian.boldt@mathematik.tu-chemnitz.de}
\author{Peter Stollmann}
\address{Fakultät f\"ur Mathematik\\Technische Universit\"at Chemnitz\\Germany}
\email{peter.stollmann@mathematik.tu-chemnitz.de }
\author{Felix Wirth}
\address{Fakultät f\"ur Mathematik\\Technische Universit\"at Chemnitz\\Germany}
\email{felix.wirth@s2019.tu-chemnitz.de}

\begin{document}

\begin{abstract} 
    We define the speed measure $\nu$ for mappings $\gamma:I\to X$ from an interval to a metric space that are locally of bounded variation. We characterize continuity and absolute continuity of $\gamma$ in terms of $\nu$ and identify the Radon-Nikod\'ym derivative of $\nu$ with respect to Lebesgue measure as the metric speed of $\gamma$. In doing so we prove an extension of the Banach-Zaretsky theorem. 
\end{abstract}


\maketitle

\section{Introduction}

In this note we define the \emph{speed measure} $\nu_\gamma$ for maps $\gamma:I\to X$, where $I\subseteq\R$ is an interval, $X$ a metric space and $\gamma$ is of locally bounded variation (see below for detailed definitions). For curves, i.e. under the additional assumption of continuity, this measure had already been studied, at least for compact $I$, in \cite{heinonen}. In the more general case we consider here, $\nu_\gamma$ measures, broadly speaking, jump sizes plus length. Therefore, curves of locally bounded variation can be characterized as those maps for which the speed measure is continuous, i.e. has no atoms.

In Section 3 we go on to characterize absolutely continuous curves as those $\gamma$ for which $\nu_\gamma$ is absolutely continuous w.r.t.\ Lebesgue measure. This can be viewed as a version in the metric space valued setting, of the Banach-Zaretsky theorem, named after the result proven independently by Banach and Zaretsky in \cite{banach} and \cite{zar}, respectively.  The original result characterizes absolute continuity for continuous real-valued functions of bounded variation in terms of Luzin's property (N), which  expresses absolute continuity of the one-dimensional Hausdorff measure of the curve w.r.t. Lebesgue measure. A generalization to curves with values in metric spaces was given by Duda and Zajicek in \cite{duda}; we discuss how our approach can be used to deduce this latter result.

For curves on compact intervals the equivalence of metric and measure theoretic absolute continuity was already noticed in \cite{heinonen}, Proposition 4.4.25, and we couldn't agree more with the authors of the latter book, who write in the paragraph preceeding Proposition 4.4.25 that this equivalence is ``a direct consequence of the definitions'' and ``basic measure theory''.

In Section 4 we introduce the metric derivative, which again was studied earlier for curves, \cite{bbi}, \cite{heinonen}, \cite{kirchheim}. We show that it is the Radon-Nikod\'ym derivative of the absolutely continuous part of the speed measure. This explains where this derivative exists: namely, the complement of the set of points where it exists is a support of the singular part of the speed measure.

We want to emphasize that the generalizations we present are natural and the proofs are all easy and straightforward due to the measure theoretic approach we follow. In this sense, we regard an advertisement for this measure theoretic approach as one of the main points of this note.

The results presented are partly taken from the Bachelor's thesis of the last named author.  

We dedicate this article to the memory of Peter Harmand, whose untimely death we deeply mourn. 

\section*{acknowledgement}
We are grateful to the anonymous referee for pointing out a gap in a previous version of our article.

\section{The set-up}

In the following $I\subseteq \R$ will denote an interval, i.e.\ a connected subset of $\R$. We write $\mathring I$ for the interior of $I$, so that $\mathring I = (a,b)$ with $a = \inf I\in\R\cup\{-\infty\}$ and $b = \sup I\in\R\cup\{+\infty\}$.

We will be concerned with mappings $\gamma:I\to X$, where $(X,\dist)$ is a metric space; as usual, $\gamma$ will be called a \emph{curve} if it is continuous. We next introduce the \emph{variation} of $\gamma$ on subintervals $J\subseteq I$:
\begin{equation*}
    \Var(\gamma;J)\coloneqq \sup\left\{\Var(\gamma;\mathcal Z)\coloneqq \sum_{k=1}^n\dist(\gamma(t_k),\gamma(t_{k+1}))\,\big|\,\mathcal Z=(t_1,\ldots,t_{n+1})\in \Par(J)  \right\}\in [0,\infty]\,,
\end{equation*}
where $\Par(J)\coloneqq \left\{(t_1,\ldots,t_{n+1})\in J^{n+1}\,|\, n\in\N, t_1\leq t_2\leq\ldots\leq t_{n+1}  \right\}$ denotes the set of partitions of $J$. We say that $\gamma$ is of \emph{bounded variation} on $J$, if $\Var(\gamma;J)<\infty$ and denote
\begin{equation*}
    \BV(I;X) \coloneqq \left\{\gamma:I\to X\,|\,\Var(\gamma;I)<\infty \right\}\,.
\end{equation*}
Furthermore, we say $\gamma$ is \emph{locally of bounded variation} if it is an element of
\begin{equation*}
    \BVl(I;X) \coloneqq \left\{\gamma:I\to X\,|\,\forall a,b\in I \text{ with } a\leq b: \Var(\gamma;[a,b])<\infty \right\}\,.  
\end{equation*}

\begin{remark}\label{basic}
    \begin{enumerate}[(1)]
        \item For curves $\gamma\in \Csf(I;X)$, it is clear that
        \[
            \Var(\gamma;J) = \operatorname{L}(\gamma_{|J})\,,
        \]
        where $\Lsf$ denotes the length of $\gamma$ and where it is customary to speak of \emph{rectifiability} instead of \emph{bounded variation}; see \cite{bbi}, Section 2.3.2 for properties of $\operatorname{L}$.
        \item In the general case, the following properties are immediate from the definition:
            \begin{enumerate}[(a)]
                \item Monotonicity: for any interval $J\subseteq I$: $\Var(\gamma;J)\leq \Var(\gamma;I)$.
                \item Additivity: for $c\in I$: $\Var(\gamma;I) = \Var(\gamma;(-\infty,c]\cap I)+\Var(\gamma;[c,\infty)\cap I)$.
                \item Inner regularity: for any interval $J\subseteq I$: $\Var(\gamma;J)=\sup\left\{\Var(\gamma;[a,b])\,|\,[a,b]\subseteq J \right\}$.
            \end{enumerate}
        \item Again in the general case, we introduce the convention $\Var(\gamma;[b,a])\coloneqq -\Var(\gamma;[a,b])$ for $a,b\in I, a\leq b$ (which is consistent since $\Var(\gamma;\{a\})=0$ for every $a\in I$). For $c\in I$ we define
        \begin{equation*}
            \V_\gamma\coloneqq \V_{\gamma,c}:I\to \R\,,\quad \V_\gamma(t)\coloneqq \Var(\gamma;[c,t])\,.
        \end{equation*}
        By the previous remark, $\V_\gamma$ is nondecreasing and
        \[
                \phantom{\qquad (a,b\in I)}  \Var(\gamma;[a,b]) = \V_\gamma(b)-\V_\gamma(a) \qquad (a,b\in I)\,.
        \]
    \end{enumerate}
\end{remark}

\begin{definition}
    For $\gamma:I\to X$ set
    \[
        \Dsf(\gamma)\coloneqq \left\{t\in I\,|\, \gamma \text{ is not continuous at }  t\right\}\,.
    \]
\end{definition}

The following properties are evident:
\begin{remark}\label{rem:evident} Let $\gamma\in \BVl(I;X)$. Then
    \begin{enumerate}[(1)]
        \item $\V_\gamma$ admits left and right limits, as usual denoted by $\V_\gamma(t-)$ and $\V_\gamma(t+)$, whenever $t\in (\inf I,\infty)\cap I$, $t\in (-\infty,\sup I)\cap I$, respectively. To ease notation, we define $\V_\gamma(a-)\coloneqq \V_\gamma(a)$ if $a=\inf I\in I$ and $\V_\gamma(b+)\coloneqq \V_\gamma(b)$ if $b=\sup I\in I$.
        \item If $t\in I$ and $(t_n)_{n\in\N}$ is a sequence in $(-\infty,t)\cap I$ or in $(t,+\infty)\cap I$ such that $t_n\to t$ as $n\to\infty$, then $(\gamma(t_n))_{n\in\N}$ is a Cauchy sequence. Hence, if $X$ is complete, then $\gamma(t-)$ and $\gamma(t+)$ exist. As above, we define $\gamma(a-)\coloneqq \gamma(a)$ if $a = \inf I\in I$ and $\gamma(b+)\coloneqq \gamma(b)$ if $b=\sup I\in I$.
    \end{enumerate}
\end{remark}

The next proposition links the left- respectively right-continuity of $\gamma$ to that of $\V_\gamma$.
\begin{proposition}\label{prop} Let $\gamma\in \BVl(I;X)$ and $t\in I$. Let $(s_n)_{n\in\N}$ be a sequence in $(-\infty,t)\cap I$ with $s_n\to t$ as $n\to\infty$. Then,
    $$
    \lim_{n\to\infty}\Var(\gamma;[s_n,t]) = \lim_{n\to\infty}\dist(\gamma(s_n),\gamma(t)) .
    $$
    If $X$ is complete, this limit agrees with $\dist(\gamma(t-),\gamma(t))$. An analogous statement holds if $(s_n)_{n\in\N}$ is a sequence in $(t,+\infty)\cap I$. Hence, $\gamma$ is left-, respectively right-continuous at $t$ $\iff$ $\V_\gamma$ is left-, respectively right-continunous at $t$. In particular,
        \[
            \Dsf(\gamma) = \Dsf(\V_\gamma)\,.
        \]

\end{proposition}
\begin{proof}
Pick partitions $(\tau_1,\ldots,\tau_{m_n+1})$ of $[s_n,t]$ so that
$$
\Var(\gamma;[s_n,t])-\sum_{k=1}^{m_n}\dist(\gamma(\tau_k),\gamma(\tau_{k+1}))\le\frac{1}{n}\mbox{  for  }n\in\N .
$$
Without restriction, $\tau_{m_n}<\tau_{m_n+1}=t$; put $t_n:=\tau_{m_n}$. Then
\begin{align*}
    \left|\Var(\gamma;[s_n,t])-\dist(\gamma(s_n),\gamma(t))\right| &\le 
    \frac{1}{n}+\left|\sum_{k=1}^{m_n}\dist(\gamma(\tau_k),\gamma(\tau_{k+1}))-\dist(\gamma(t_n),\gamma(t))\right| + \dist(\gamma(s_n),\gamma(t_n))\\
    &\le \frac{1}{n}+\sum_{k=1}^{m_n-1}\dist(\gamma(\tau_k),\gamma(\tau_{k+1}))+ \dist(\gamma(s_n),\gamma(t_n))\\
    &\le \frac{1}{n}+\Var(\gamma;[s_n,t_n])+ \dist(\gamma(s_n),\gamma(t_n))\,.
\end{align*}
Since $t>t_n\geq s_n\to t$ as $n\to\infty$, we have by Remark~\ref{basic}(3) and Remark~\ref{rem:evident}(1),
\[
    \dist(\gamma(s_n),\gamma(t_n))\leq \Var(\gamma;[s_n,t_n]) = \V_\gamma(t_n)-\V_\gamma(s_n)\to \V_\gamma(t-)-V_\gamma(t-) = 0\,,
\]
which completes the proof.
\end{proof}
We shall now point out why $\Var(\gamma;\cdot)$ is, in general, not additive: To ease notation, we write, for $\gamma\in\BVl(I;X)$, $t\in \mathring I$,
\begin{align*}
\dist(\gamma(t-),\gamma(t))&:=\lim_{s\nearrow t}\dist(\gamma(s),\gamma(t)),\\
    \dist(\gamma(t),\gamma(t+))&:=\lim_{s\searrow t}\dist(\gamma(t),\gamma(s)) .
\end{align*}
These limits exist by the previous proposition; by passing to the completion of $X$ if necessary, '$:=$' is in fact an equality.

\begin{corollary} Let $\gamma\in \BVl(I;X)$.
    \begin{itemize}
        \item[(1)] If $t\in\Dsf(\gamma)$ and $a\in I$, $a<t$, and $\gamma$ is not left-continuous at $t$, then
                \begin{align*}
                    \Var(\gamma;(a,t]) &= \Var(\gamma;(a,t)) + \dist(\gamma(t-),\gamma(t))\\
                        & > \Var(\gamma;(a,t))+\Var(\gamma;\{t\})\,,
                \end{align*}
                as $\Var(\gamma;\{t\})=0$ by definition.
        \item[(2)] If $t\in\Dsf(\gamma)$ and $b\in I$, $t < b$, and $\gamma$ is not right-continuous at $t$, then
                \begin{align*}
                    \Var(\gamma;[t,b)) &= \Var(\gamma;(t,b)) + \dist(\gamma(t),\gamma(t+))\\
                        & > \Var(\gamma;(t,b))+\Var(\gamma;\{t\})\,.
                \end{align*}
        \end{itemize}
\end{corollary}
\begin{proof}
    (1)  Pick a sequence $(t_n)_{n\in\N}$ in $(a,t)$ that increases to $t$. By additivity and inner regularity of $\Var$, see Remark \ref{basic},
    \begin{align*}
        \Var(\gamma;(a,t])&=\lim_{n\to\infty}\left(\Var(\gamma;(a,t_n])+\Var(\gamma;[t_n,t])\right)\\
        &=\Var(\gamma;(a,t))+\lim_{n\to\infty}\Var(\gamma;[t_n,t])\\
        &=\Var(\gamma;(a,t)) + \dist(\gamma(t-),\gamma(t)),
    \end{align*}
    where we used the previous proposition in the last step. Assertion (2) is proven in the same way.
\end{proof}
We are now ready to introduce the main character of this note.

\begin{definition} Let $\gamma\in \BVl(I;X)$. 
    We denote by $\vv(t)\coloneqq\vv_{\gamma}(t)\coloneqq \V_{\gamma}(t+)$, $t\in I$, the right-continuous modification of $\V_{\gamma}$ and by $\mu_{\vv}$ the corresponding  Lebesgue-Stieltjes measure  on $(I,\mathcal B)$. In case $a=\inf I\in I$, we define  $\nu:=\mu_{\vv}+ (\V_{\gamma}(a+)-\V_{\gamma}(a))\delta_a$ and else $\nu:=\mu_{\vv}$ and call $\nu$ the \emph{speed measure} of $\gamma$.
\end{definition}

The existence of such a Lebesgue-Stieltjes measure $\mu_{\vv}$ for any right-continuous, non-decreasing function $\vv$  is a classical fact; it is characterized by
    \begin{equation}\label{lsmeasure}
        \phantom{\qquad (s,t\in I, s\leq t)}        \mu_{\vv}((s,t])  = \vv(t)-\vv(s) \qquad (s,t\in I, s\leq t)\, 
    \end{equation} 
    and can be constructed by the Carath\'eodory extension theorem, by prescribing the values given in \eqref{lsmeasure} on the semiring of semiopen intervals, see e.g. \cite{bc}, Theorem~3.5.1 or \cite{klenke}, Example 1.56. There it is stated for functions $\vv$ defined on $\R$; the slight generalization to the case at hand presents no difficulties. Observe that $\vv_\gamma$ can be regarded as a shifted distribution function of $\nu_\gamma$. In a classical probabilistic context, one would and could require $I=\R$, $\vv_\gamma(t)\to 0 $ as $t\to-\infty$, $\vv_\gamma(t)\to+ 1$ as $t\to+\infty$, see \cite{bc}, p.\ 356. Of course, this is not possible for the general set-up here.

\begin{theorem}\label{thm:speed-measure-extends-variation}
    Let $\gamma\in\BVl(I;X)$. Then the speed measure $\nu$ of $\gamma$ is the unique measure on $(I,\mathcal B)$ such that for all intervals $J\subseteq I$, $J$ open in $I$, one has
    \begin{equation}\label{eqn:speed-measure-var}
        \nu(J) = \Var(\gamma;J)\,.
    \end{equation}
    Moreover, for $t\in I$
    \begin{gather}\label{eqn:speed-measure-points}
    \begin{aligned}
        \nu(\{t\}) & = \V_\gamma(t+)-\V_\gamma(t-) \\
        & = \dist(\gamma(t-),\gamma(t)) + \dist(\gamma(t),\gamma(t+))   \,.
    \end{aligned}
    \end{gather}
    Furthermore, $\nu$ is continuous if and only if $\gamma$ is a curve. In this case, \eqref{eqn:speed-measure-var} holds for all intervals $J\subseteq I$.
\end{theorem}
\begin{proof}
    Denote $\mathscr G = \left\{(s,t]\,|\,s,t\in I \right\}\cup \{\{a\}\}$ if $a=\inf I\in I$ and $\mathscr G = \left\{(s,t]\,|\,s,t\in I \right\}$ else. Then $\mathscr G$ is a semiring which generates $\mathcal B(I)$. Every measure that satisfies \eqref{eqn:speed-measure-var} is $\sigma$-finite and uniquely determined on $\mathscr G$. Therefore, the uniqueness results from Theorem 1.53 in \cite{klenke}.

    The first part of \eqref{eqn:speed-measure-points} follows by definition of $\nu$ through $\vv_\gamma$, the second equality follows from Proposition \ref{prop} above.

    By \eqref{eqn:speed-measure-points}, it is clear that $\gamma$ is continuous iff $\nu$ is continuous.
\end{proof}

\begin{remark}\label{rem:meas-var}
    \begin{enumerate}[(1)]
        \item The sum in \eqref{eqn:speed-measure-points} gives the total jump size of $\gamma$ at $t$. Keeping this in mind,  \eqref{eqn:speed-measure-var} above immediately implies that $\Var(\gamma;J)\leq \nu(J)$ for all intervals $J\subseteq I$, with equality iff $J\setminus\mathring J\subseteq I\setminus\Dsf(\gamma)$.
        \item The speed measure for curves already appeared in \cite{heinonen}, Theorem~4.4.8.
    \end{enumerate}
\end{remark}

\section{Absolute continuity and the Banach-Zaretsky Theorem}

Before stating the main result of this section, we need to introduce some notions even though they are probably well-known to most readers. 

\begin{definition}
    \begin{enumerate}[(1)]
        \item Absolute continuity of measures: let $\mu,\nu$ be measures on the same measurable space $(\Omega,\mathcal F)$. Then $\nu$ is absolutely continuous w.r.t.\ $\mu$, denoted $\nu\ll\mu$, if: $\forall F\in\mathcal F$: $\mu(F)=0 \implies \nu(F)=0$.
        \item Absolute continuity of mappings $\gamma:I\to X$, where $I\subseteq \R$ is an interval and $(X,\dist)$ a metric space: $\gamma$ is absolute continuous if $\forall\mkern6mu \varepsilon>0\mkern6mu\exists\mkern6mu\delta > 0\mkern6mu\forall\mkern3mu (a_k,b_k), k=1,\ldots n$ pairwise disjoint subintervals of $I$  
        \begin{equation}\label{eqn:absolute-continuity-maps}
            \sum_{k=1}^n (b_k-a_k) \leq \delta\implies \sum_{k=1}^n\dist(\gamma(b_k),\gamma(a_k))\leq \varepsilon\,. 
        \end{equation}
        We write 
        $$\AC(I;X) \coloneqq \left\{ \gamma:I\to X\,|\, \gamma \text{ is absolutely continuous} \right\}$$ 
        
        and 
        $$\ACl(I;X)\coloneqq\left\{\gamma:I\to X\,|\, \gamma_{|[a,b]}\in \AC([a,b];X) \text{ for all } a,b\in I, a\leq b \right\} .$$
    \end{enumerate}
\end{definition}

We are now ready to state

\begin{theorem}[Banach-Zaretsky Theorem - speed measure version]\label{thm:Banach-Zaretsky}
    Let $\gamma:I\to X$. Then $\gamma\in\ACl(I;X)$ if and only if $\gamma \in \BVl(I;X)$ and its speed measure $\nu$ is absolutely continuous with respect to Lebesgue measure $\lambda$.
\end{theorem}

Given the measure theoretic framework, the proof is rather elementary. We single out one crucial observation that links metric absolute continuity to measure theoretic absolute continuity.

\begin{lemma}
    Let $(\Omega,\mathcal F)$ be a measure space, $\mu, \nu$ measures on $(\Omega, \mathcal F)$ with $\mu$ finite. Then the following are equivalent:
    \begin{enumerate}[(i)]
        \item $\nu \ll\mu$.
        \item $\forall\mkern3mu \varepsilon > 0\mkern6mu\exists\mkern6mu\delta > 0\mkern6mu\forall\mkern3mu F\in\mathcal F: \mu(F)\leq \delta \implies\nu(F)\leq\varepsilon$.
    \end{enumerate}
\end{lemma}
\begin{proof}
  The implication ``$(ii)\implies(i)$'' is trivial, and ``$(i)\implies (ii)$'' is contained in \cite{bc}, Theorem 5.2.7.
\end{proof}

Now the reason for the equivalence in the Banach-Zaretsky theorem is clearly visible: in \eqref{eqn:absolute-continuity-maps} the left-hand side can be reinterpreted as $\lambda(\cup_k[a_k,b_k])\leq \delta $. On the right-hand side we have 
\begin{equation}\label{eqn:abs-cts-rhs-meas}
    \sum_{k=1}^n\dist(\gamma(a_k),\gamma(b_k)) \leq \sum_{k=1}^n \Lsf\left(\gamma_{|[a_k,b_k]}\right) = \sum_{k=1}^n\nu([a_k,b_k])\,. 
\end{equation}

\begin{proof}[Proof of Theorem 3.2]
    $[\gamma\in\BVl(I;X)\wedge \nu\ll\lambda]\implies \gamma\in\ACl(I;X)$ immediately follows from the above observation: for $J$ a compact interval in $I$, \eqref{eqn:abs-cts-rhs-meas} gives, in light of the above lemma, absolute continuity of $\gamma$ on $J$.\\
    
    \noindent
    $\gamma\in\ACl(I;X)\implies[\gamma\in\BVl(I;X)\wedge \nu\ll\lambda]$: Let $J\subseteq I$ be a compact interval. The definition readily implies $\AC(J;X)\subseteq \BV(J;X)\cap \Csf(J;X)$, so that we are left to prove $\nu\ll\lambda$. Since we can write $I=\bigcup_{m\in\N} J_m$, $J_m$ compact, it suffices to consider $I$ compact. 

    We will first show that for every $\varepsilon>0$ there exists $\delta>0$ so that for $E=\bigcup_{k=1}^n[a_k,b_k]\subseteq I$ with $\lambda(E)\leq \delta$ we have $\nu(E)\leq \varepsilon$. Let $\varepsilon > 0$, pick $\delta>0$ as in the definition of absolute continuity of $\gamma$ and let $E=\bigcup_{k=1}^n[a_k,b_k]$ with pairwise disjoint $(a_k,b_k)$ and $\lambda(E)\leq\delta$. For each $k=1,\ldots,n$, pick $\mathcal Z_k=\left(t^k_1,\ldots,t^k_{m(k)+1}\right)\in\Par([a_k,b_k])$. Since
    \[
        \sum_{k=1}^n\sum_{\ell=1}^{m(k)}|t_{\ell+1}-t_\ell| \leq \sum_{k=1}^n(b_k-a_k)\leq \delta\,,
    \]
    we obtain
    \[
        \sum_{k=1}^n\Lsf(\gamma,\mathcal Z_k) = \sum_{k=1}^n\sum_{\ell=1}^{m(k)}\dist(\gamma(t^k_{\ell+1}), \gamma(t^k_\ell))\leq \varepsilon\,.
    \]
    Taking the supremum over all $\mathcal Z_k\in\Par([a_k,b_k])$, $k=1,\ldots,n$, we thus get
    \[
        \nu(E) = \sum_{k=1}^n\Lsf\left(\gamma_{|[a_k,b_k]}\right) \leq \varepsilon\,.
    \]
    To finish the proof that $\nu\ll\lambda$, let $N\in\mathcal B(I)$ with $\lambda(N)=0$. By the regularity of $\lambda$ there exist $a_i,b_i\in I$ such that $N\subseteq \cup_{i=1}^\infty[a_i,b_i]$ and $\sum_{i=1}^\infty(b_i-a_i)\leq\delta$. Then we have by the above
    \[
        \nu(N) \leq \nu\left(\bigcup_{i=1}^\infty [a_i,b_i]\right) = \lim_{n\to\infty} \nu\left(\bigcup_{i=1}^n[a_i,b_i] \right)\leq \varepsilon\,.
    \]
\end{proof}

A curve $\gamma\in \BVl(I;X)$ is said to have \emph{Luzin's property (N)} if 
\begin{equation*}
    \forall N\subseteq I: \lambda(N) = 0 \implies \mathscr H^1(\gamma(N))=0\,,
\end{equation*}
where $\mathscr H^1$ denotes the 1-dimensional outer Hausdorff measure of $(X,\dist)$. Given the elementary inequality for curves, see Theorem 2.6.2 and Remark 2.6.3 in \cite{bbi}:
\[
    \mathscr H^1(\gamma(A))\leq \nu(A)
\]
for all $A\in\mathcal B(I)$, with equality if $\gamma$ is a simple curve, we obtain the following corollary to Theorem~\ref{thm:Banach-Zaretsky}.
\begin{corollary}
    Let $\gamma:I\to X$. If $\gamma\in \ACl(I;X)$, then $\gamma\in \Csf(I;X)\cap \BVl(I;X)$ and $\gamma$ has Luzin's property (N). If $\gamma$ is simple, the converse holds.
\end{corollary}
This gives a particularly simple proof of the Banach-Zaretsky theorem, which we now recall for convenience of the reader, for simple curves. In the present generality it is due to Duda and Zajicek, \cite{duda} to which we refer for further references on earlier extensions.

\begin{theorem}[Banach-Zaretsky Theorem]\label{thm:B-Z}
  Let $\gamma:I\to X$. Then $\gamma\in\ACl(I;X)$ if and only if $\gamma \in C(I;X)\cap \BVl(I;X)$ and it satisfies Luzin's property (N).
\end{theorem}
Here is how this Theorem can be obtained quite quickly from our version above:  Using the fact that length is absolutely continuous w.r.t.\ the Hausdorff measure of the image, see e.g., Exercise 2.6.4 in \cite{bbi},  Theorem 2.10.13 in \cite{federer}, or the corresponding discussion in \cite{duda07}, a curve with Luzin's property (N) has absolutely continuous speed measure with respect to Lebesgue measure and, consequently is locally absolutely continuous. The converse direction was adressed above. 

\section{The Lebesgue Decomposition Theorem and the metric derivative}

The following definition seems to have appeared first in \cite{kirchheim} in the context of Lipschitz maps $\R^n\to X$.

\begin{definition}
    For $\gamma\in \BVl(I;X)$ we define the \emph{metric derivative of $\gamma$ at $t\in I$} as
    \begin{equation*}
        |\dot\gamma|(t)\coloneqq \lim_{\varepsilon\to 0}\frac{\dist(\gamma(t+\varepsilon),\gamma(t))}{|\varepsilon|}
    \end{equation*}
    provided the limit exists.
\end{definition}
We will show below, that the metric derivative exists $\lambda$-a.e.\ and is the Radon-Nikod\'ym derivative of the absolutely continuous part of the speed measure of $\gamma$. To do so, we need the following generalization of a result known for rectifiable curves, see \cite{bbi}, Theorem\ 2.7.4. In a previous version of this paper, we used this result in the present more general setting of mappings of locally bounded variation, not taking into account that the proof given in \cite{bbi} uses continuity in a crucial way. We thank the referee for pointing out this oversight and  mention that for curves, our proof below simplifies considerably.
\begin{theorem}
Let $I=[a,b]$ and $\gamma\in\BV(I;X)$. Then, for Lebesgue almost all $t\in [a, b]$
\begin{equation}\label{eqn:either-or}
     \qquad \liminf_{\delta,\delta'\searrow 0} \frac{\Var\left(\gamma; [t-\delta,t+\delta']\right)}{\delta+\delta'} = 0 \qquad \text{ or } \qquad \lim_{\delta,\delta' \searrow 0} \frac{\dist(\gamma(t-\delta),\gamma(t+\delta'))}{\Var\left(\gamma; [t-\delta,t+\delta']\right)} = 1\,. 
\end{equation}
\end{theorem}
\begin{proof}
Let $\alpha > 0$ and consider the set $A_\alpha$ of all $t\in I$ such that
   \begin{equation*}
       \liminf_{\delta,\delta'\searrow 0} \frac{\Var\left(\gamma; [t-\delta,t+\delta']\right)}{\delta+\delta'} > \alpha \quad \text{ and } \quad \lim_{\delta,\delta' \searrow 0} \frac{\dist(\gamma(t-\delta),\gamma(t+\delta'))}{\Var\left(\gamma; [t-\delta,t+\delta']\right)} < 1 - \alpha .
   \end{equation*}
Since the alternative claimed in the theorem holds for all $t\not\in \cup_{n\in\N}A_{\frac{1}{n}}$, it suffices to prove that $\lambda(A_\alpha)=0$ for every $\alpha>0$. So fix $\alpha>0$ and $A=A_\alpha$.

(1) Fix $\varepsilon>0$ and let $S=\{s_1, ... , s_L\}\subset I$ be finite with $s_1=a, s_L=b$ in increasing order so that
\begin{equation}\label{eq0}
        \sum_{t\not\in S}\nu(\{t\}) < \frac{\varepsilon\alpha^2}{8}\,.     
\end{equation}
Such an $S$ exists, since $\nu$ is finite and so 
$$
\sum_{t\in I}\nu(\{t\})<\infty\, .
$$
Regularity of the Lebesgue measure gives that we find compact intervals $I_\ell\subset (s_\ell,s_{\ell+1})$ so that
\begin{equation}\label{eq1}
      \lambda(A)\le 2\sum_{\ell=1}^L \lambda(A\cap I_\ell)\,.
\end{equation}
(2) For $\ell=1, ..., L-1$ we now pick a partition $\mathcal Z\in\Par(I_\ell)$  such that 
   \begin{equation}\label{eq2}
       \Var(\gamma; I_\ell) - \Var(\gamma; \mathcal Z) < \frac{\varepsilon\alpha^2}{8}2^{-\ell}\,.
   \end{equation} 
Denote by $Z$ the set of points belonging to this partition and define
   \begin{multline*}
       \mathfrak B \coloneqq \{ [t-\delta,t+\delta']\,|\, t\in A\cap I_\ell,[t-\delta,t+\delta']\cap Z=\emptyset,\\\Var(\gamma;[t-\delta,t+\delta'])>\alpha(\delta+\delta'),\dist(\gamma(t-\delta),\gamma(t+\delta'))<(1-\alpha)\Var(\gamma;[t-\delta,t+\delta']) 
       \}\,.
   \end{multline*}
   By definition of $A$ and $\mathfrak B$, every $t\in (A\cap I_\ell)\setminus Z$ is contained in an arbitrarily short element of $\mathfrak B$, so that Vitali's covering theorem \cite{bbi}, Theorem 1.7.14 yields, a countable collection $\{[t_i-\delta_i,t_i+\delta_i']\,|\,i\in\Lambda\}$ of pairwise disjoint intervals that cover $(A\cap I_\ell)\setminus Z$ up to a Lebesgue null set. In particular,
   \[
    \sum_{i\in\Lambda}(\delta_i+\delta'_i) = \lambda\left(\bigcup_{i\in\Lambda}[t_i-\delta_i,t_i+\delta_i'] \right) \geq \lambda((A\cap I_\ell)\setminus Z)=\lambda(A\cap I_\ell)\,.
   \]
   Hence, there exist $i_1,\ldots,i_M\in\Lambda$ such that
   \begin{equation}\label{eqn:sum-mu-half}
        \sum_{j=1}^M \left(\delta_{i_j}+\delta'_{i_j}\right) \ge \frac12 \lambda(A\cap I_\ell)\,.
   \end{equation}
  (3) Since the intervals $[t_{i_j}-\delta_{i_j},t_{i_j}+\delta'_{i_j}]$, $j=1,\ldots,M$, are pairwise disjoint and do not have points of $Z$ in their interior, we can include their endpoints in a partition $\mathcal Y=(y_1,\ldots,y_N)\in\Par(I)$ refining $\mathcal Z$ by possibly adding points not lying in any of the given intervals and the remaining ones of $Z$.
   Observe that the definition of $A$ and $\mathfrak B$ gives that for every $j$ as above
   \begin{equation}\label{eq3}
       \Var(\gamma;[t_{i_j}-\delta_{i_j},t_{i_j}+\delta'_{i_j}]) -\dist(\gamma(t_{i_j}-\delta_{i_j}),\gamma(t_{i_j}+\delta'_{i_j}))\ge\alpha^2\left(\delta_{i_j}+\delta'_{i_j}\right) .
   \end{equation}
(4) We can now estimate $\lambda(A\cap I_ \ell)$ using the previous inequality: to ease notation, let us call the intervals $J_1:=[y_1,y_2]$, $J_i:=(y_i,y_{i+1}]$, for $I=1, ... , N-1$ and note that they form a disjoint decomposition of $I_\ell$: Moreover, we let $G$ be the set of those
$i$, for which $y_i=t_{i_j}-\delta_{i_j}$ for some $j$ as above. We get
\begin{align*}
    \lambda(A\cap I_\ell) &\le 2 \sum_{j=1}^M \left(\delta_{i_j}+\delta'_{i_j}\right)\\
    &\le \frac{2}{\alpha^2}\sum_{j=1}^M \left(\Var(\gamma;[t_{i_j}-\delta_{i_j},t_{i_j}+\delta'_{i_j}]) -\dist(\gamma(t_{i_j}-\delta_{i_j}),\gamma(t_{i_j}+\delta'_{i_j}))\right)\\
    &\le \frac{2}{\alpha^2}\left[ \sum_{i\in G} \left(\Var(\gamma;\overline{J_i}) -\dist(\gamma(y_i);\gamma(y_{i+1}))\right)  \right]\\
    &\le  \frac{2}{\alpha^2}\left[ \sum_{i=1}^{N-1} \left(\Var(\gamma;\overline{J_i}) -\dist(\gamma(y_i);\gamma(y_{i+1}))\right) \right] \,.
\end{align*}
Note that for the last inequality it was essential to have the closure $\overline{J_i}$ in order to assure that all terms in the sum are nonnegative. This extra complication disappears for curves, for which simply $\Var(\gamma;\overline{J_i})=\Var(\gamma;J_i)$.  However, the difference of theses terms, which is just $\nu(\{y_i\})$ for $i\ge 2$ and zero for $i=1$ is controlled by the point part of $\nu$ and we can estimate further:
\begin{align*}
    \lambda(A\cap I_\ell) 
    &\le  \frac{2}{\alpha^2}\left[ \sum_{i=1}^{N-1} \left(\Var(\gamma;J_i) -\dist(\gamma(y_i);\gamma(y_{i+1}))\right) +\sum_{t\in I_\ell}\nu(\{ t\})\right] \\
    &= \frac{2}{\alpha^2}\left[\Var(\gamma;I_\ell)-\Var(\gamma;\mathcal Y)+\sum_{t\in I_\ell}\nu(\{ t\})\right] \\
    &\le \frac{2}{\alpha^2}\left[\frac{\varepsilon\alpha^2}{8}2^{-\ell} +\sum_{t\in I_\ell}\nu(\{ t\})\right] 
\end{align*}
in view of \eqref{eq2} above. We can now finish the proof by using inequalities \eqref{eq0} and
\eqref{eq1} above:
\begin{align*}
    \lambda(A) &\le 2\sum_{\ell=1}^L \frac{2}{\alpha^2}\left[\frac{\varepsilon\alpha^2}{8}2^{-\ell} +\sum_{t\in I_\ell}\nu(\{ t\})\right] \\
    &\le \varepsilon \, .
\end{align*}
\end{proof}

\begin{corollary}\label{cor:curves-infini-geod}
Let $I$ be an arbitrary interval and $\gamma\in\BVl(I;X)$. Then, for Lebesgue almost all $t\in I$
\begin{equation*}
    \liminf_{\delta\to 0} \frac{\Var\left(\gamma;[t,t+\delta]\right)}{\delta} = 0 \qquad \text{ or } \qquad \lim_{\delta\to 0} \frac{\dist(\gamma(t),\gamma(t+\delta))}{\left|\Var\left(\gamma;[t,t+\delta]\right)\right|} = 1\,. 
\end{equation*}
\end{corollary}
\begin{proof}
    The statement follows from the previous theorem by, e.g., considering a countable exhaustion of $I$ by compact intervals.
\end{proof}

\begin{theorem}\label{thm:lebesgue-decomp}
    Let $\gamma\in\BVl(I;X)$, $\nu$ the speed measure of $\gamma$ and $\nu = \nu_\text{ac}+\nu_\text{sing}$ its Lebesgue decomposition with respect to Lebesgue measure $\lambda$. Then $|\dot\gamma|(t)$ exists for $\lambda$-almost all $t\in I$, $\nu_\text{ac}=|\dot\gamma|\cdot \lambda$, and for all intervals $J\subseteq I$ one has
    \begin{equation*}
        \Var(\gamma; J)\leq\int_J\!|\dot\gamma|\,\dd\lambda + \nu_\text{sing}(J)\,,
    \end{equation*}
    with equality if $J\setminus \mathring J\subseteq I\setminus \Dsf(\gamma)$. In particular, if $\gamma$ is a curve, one has for all intervals $J\subseteq I$:
    \begin{equation*}
        \Lsf\left(\gamma_{|J}\right) = \int_J\!|\dot\gamma|\,\dd\lambda + \nu_{\text{sing}}(J)\,.
    \end{equation*}
\end{theorem}
\begin{proof}
    In light of Remark~\ref{rem:meas-var}, we only have to show that $|\dot\gamma|(t)$ exists for $\lambda$-almost all $t\in I$ and that $|\dot\gamma|$ is the Radon-Nikod\'ym derivative of $\nu_\text{ac}$.

    To that end, note that $\vv$ is a monotonically nondecreasing function, hence differentiable $\lambda$-almost everywhere. Taking into account that $\vv$ is not differentiable at any $t\in \Dsf(\gamma)=\Dsf(\vv)$, we obtain for any $t\in I$ at which $\vv$ is differentiable
    \begin{align*}
        \vv'(t) = \lim_{\varepsilon\to 0} \frac{\vv(t+\varepsilon)-\vv(t)}{\varepsilon} = \lim_{\varepsilon\to 0} \frac{\V_\gamma((t+\varepsilon)+)-\V_\gamma(t)}{\varepsilon}\,.
    \end{align*}
    Since for any null sequence $(\varepsilon_n)_n$
    \begin{multline*}
        \lim_{n\to\infty} \frac{\V_\gamma(t+\varepsilon_n)-\V_\gamma(t)}{\varepsilon_n}\leq \lim_{n\to\infty} \frac{\V_\gamma((t+\varepsilon_n)+)-\V_\gamma(t)}{\varepsilon_n}\\
        \leq \lim_{n\to\infty} \frac{\V_\gamma(t+\varepsilon_n+\varepsilon_n^2)-\V_\gamma(t)}{\varepsilon_n} = \lim_{n\to\infty} \frac{\V_\gamma(t+\varepsilon_n+\varepsilon_n^2)-\V_\gamma(t)}{\varepsilon_n (1+\varepsilon_n)}\,,
    \end{multline*}
    we have
    \begin{equation}\label{eqn:v-derivative}
        \vv'(t) =  \lim_{\varepsilon\to 0} \frac{\Var(\gamma;[t,t+\varepsilon])}{\varepsilon}=\lim_{\varepsilon\to 0} \frac{|\Var(\gamma;[t,t+\varepsilon])|}{\dist(\gamma(t+\varepsilon),\gamma(t))}\cdot\frac{\dist(\gamma(t+\varepsilon),\gamma(t))}{|\varepsilon|}\,.
    \end{equation}
    By Corollary\ \ref{cor:curves-infini-geod} we have for $\lambda$-almost all $t\in I$
    \begin{equation*}
        \liminf_{\varepsilon\to 0}\frac{\Var(\gamma;[t,t+\varepsilon])}{\varepsilon} = 0 \quad \text{ or } \quad \lim_{\varepsilon\to 0}\frac{\dist(\gamma(t+\varepsilon),\gamma(t))}{|\Var(\gamma;[t,t+\varepsilon])|} = 1\,.
    \end{equation*}
    In the first case, we have $\vv'(t) = 0$, and since $\dist(\gamma(t+\varepsilon),\gamma(t))\leq |\Var(\gamma;[t,t+\varepsilon])|$, we get $\vv'(t)=|\dot\gamma|(t)$. In the second case, the first factor in the product in \eqref{eqn:v-derivative} is one in limit, so that in this case we have $\vv'(t) = |\dot\gamma|(t)$ too.

    Next, let $f\in \Lsf^1_\text{loc}(I,\lambda)$ be the Radon-Nikod\'ym derivative of $\nu_\text{ac}$ with respect to $\lambda$. Define for $n\in\N$ and $t\in \mathring I$ the set $E_n(t)\coloneqq (t,t+\sfrac{1}{n}]$. Then $(E_n(t))_{n\in\N}$ shrinks to $t$ nicely in the sense of \cite{rudin}, 7.9. By \cite{rudin}, Theorem~7.14 we have for $\lambda$-almost all $t\in I$
    \begin{equation*}
        f(t) = \lim_{n\to\infty}\frac{\nu(E_n(t))}{\lambda(E_n(t))} = \lim_{n\to\infty} \frac{\vv(t+\sfrac{1}{n})-\vv(t)}{\sfrac 1n}=\vv'(t) = |\dot\gamma|(t)\,.
    \end{equation*}
\end{proof}

We end by showing that our result above gives a nice characterization of the classes $AC^p$ of curves introduced in \cite{ambrosio}. While there, intervals are assumed to be open, let us stick to our slightly more general set-up and define for $1\leq p\leq \infty$
$$
AC^p(I;X):=\left\{\gamma:I\to X\mid \exists g\in \Lsf^p(I)\quad\forall s,t\in I, s\le t: \dist(\gamma(s),\gamma(t))\le \int_s^t\! g\,\dd \lambda\right\}\,,
$$
which we caution the reader not to confuse with our earlier defined space $\AC(I;X)$.

First note that for $g$ as in the definition of $AC^p$, one has
$$
\forall s,t\in I, s\le t: \Lsf(\gamma|_{[s,t]}) \le \int_s^t g\,\dd\lambda\,,
$$
as is readily seen from the definition of $\mathsf L$, which is equivalent to $\nu_\gamma\leq g\cdot\lambda$. Hence, it follows from Theorem~\ref{thm:lebesgue-decomp} that $|\dot\gamma|\leq g$ $\lambda$-a.e. for any such $g$. In combination with Theorem~\ref{thm:Banach-Zaretsky} this yields $AC^p(I;X)\subseteq \ACl(I;X)$, so that we obtain:
\begin{corollary}
    $AC^p(I;X)=\{\gamma\in \ACl(I;X)\mid  |\dot\gamma|\in \Lsf^p(I)\}$, in particular $AC^1(I;X)=\ACl(I;X)\cap \BV(I;X)$. 
\end{corollary}

\end{document}